\documentclass[12pt]{article}

\usepackage{amsmath}
\usepackage{latexsym}
\usepackage{amssymb,enumerate}
\usepackage{amsfonts}

\newcommand{\proof}{\par\noindent{\it Proof.\ \ }}

\newcommand\qed{\ifmmode\square\else\nolinebreak\hfill
$\Box$\fi\par\vskip12pt}

\newcommand\calD{{\mathcal D}} 
\newcommand\calM{{\mathcal M}} 
\newcommand\calP{{\mathcal P}}
\newcommand\bfK{{\bf K}} 
 
\DeclareMathOperator\diam{{\sf diam}}
\DeclareMathOperator\dist{{\sf d}} 
\DeclareMathOperator\Aut{{\sf Aut}}
\DeclareMathOperator\GL{{\sf GL}}
\DeclareMathOperator\AGL{{\sf AGL}}

\newtheorem{theorem}{Theorem}[section]%
\newtheorem{lemma}[theorem]{Lemma}%
\newtheorem{proposition}[theorem]{Proposition}%
\newtheorem{definition}[theorem]{Definition}%
\newtheorem{problem}[theorem]{Problem}%
\newtheorem{example}[theorem]{Example}%
\newtheorem{construction}[theorem]{Construction}%
\newtheorem{remark}[theorem]{Remark}%
 \newtheorem{hypothesis}[theorem]{Hypothesis}

\title{Locally $s$-distance transitive graphs and pairwise transitive designs
\footnote{The research for this paper was supported by Australian Research Council Discovery Grant DP0770915 and Federation Fellowship Grant FF0776186. Cheryl Praeger is also affiliated with King Abulazziz University, Jeddah, Saudi Arabia.
}}
\author{Alice Devillers,
 Michael Giudici, Cai Heng Li and Cheryl E.~Praeger\\
Centre for the Mathematics of Symmetry and Computation\\
 School of Mathematics and Statistics\\
 The University of Western Australia\\
 35 Stirling Highway, Perth WA 6009, Australia\\
emails: \texttt{ alice.devillers@uwa.edu.au, michael.giudici@uwa.edu.au,} \\
 \texttt{cai.heng.li@uwa.edu.au,  cheryl.praeger@uwa.edu.au}}

 \date{}

\begin{document}
\maketitle

\begin{abstract}
The study of locally $s$-distance transitive graphs initiated by the authors in previous work, identified that graphs with a star quotient are of particular interest. 
This paper shows that the study of locally $s$-distance transitive graphs with a star quotient is equivalent to the study of a particular family of designs with strong symmetry properties that we call nicely affine and pairwise transitive. We show that a group acting regularly on the points of  such a design must be abelian and give general construction for this case.
\end{abstract}

\section{Introduction}

In \cite{DGLP}, we studied finite locally $s$-distance transitive graphs and found that graphs with a star quotient (which we call starlike here, see Definition~\ref{def:starlike}) were of particular interest. All graphs considered in this paper are finite, simple and without loops.
The parameter $s$ is a positive integer, and a graph is said to be \emph{locally $(G,s)$-distance transitive} if the graph has diameter at least $s$, admits $G$ as a subgroup of automorphisms, and for each vertex $v$ and each positive integer $i\leq s$, the stabiliser $G_v$ acts transitively on the set of vertices at distance $i$ from $v$. 
Whenever $G$ is the full automorphism group of the graph, we sometimes simply say the graph is \emph{locally $s$-distance transitive}. 

In this paper we show that for $s\geq 4$, the study of locally $s$-distance transitive graphs with a star quotient 
can be transformed into the study of designs with some specified structural and symmetry properties  (see Section 
\ref{Def} for the definition of a design). Such graphs are bipartite, and for a bipartite graph $\Gamma$ with ordered 
bipartition $(B|B')$, we define the {\it adjacency design}  ${\cal D}(\Gamma)$  to have  point set $B$ and block set $B'$, 
such that a point and a block are incident if the corresponding vertices are adjacent in $\Gamma$.
We write  ${\cal D}(\Gamma)=(B,B',E\Gamma)$ where $E\Gamma$, the edge set of $\Gamma$, is seen as a subset of $B\times B'$ (we acknowledge a little abuse of notation here).
  Conversely, given a design ${\cal D}=({\cal P}, {\cal B}, {\cal I})$, 
we define the {\it incidence graph} of $\calD$ as the  bipartite graph $\Gamma({\cal D})$ with ordered bipartition $({\cal P}| {\cal B})$ and adjacency given by incidence. So ${\cal D}(\Gamma({\cal D}))={\cal D}$ and $\Gamma({\cal D}(\Gamma))=\Gamma$.
The adjacency designs for locally $4$-distance transitive graphs turn out to be nicely affine and pairwise transitive (see definitions below).

\begin{definition}\label{def:nice} 
{\rm Let ${\cal D}=({\cal P}, {\cal B}, {\cal I})$ be a design and $N$ be an automorphism group of $\calD$. Then  $\calD$ is 
called {\it $N$-nicely affine} if $N$ is transitive on ${\cal P}$ 
and there is a constant $\mu$ such that distinct blocks  are incident with exactly $\mu$ common points if they are in different $N$-orbits and are disjoint if they are in the same $N$-orbit. 
}
\end{definition}

Since $N$ is transitive on ${\cal P}$ it follows from Definition~\ref{def:nice} that each $N$-orbit in  ${\cal B}$ is 
a \emph{parallel class}, that is, the subsets of points incident with blocks in the $N$-orbit form a partition of 
$\cal P$. If all blocks of $\cal D$ are incident with the same number of points then
an $N$-nicely affine design is affine in the usual sense (see Section~\ref{designs}), 
but we see in Example~\ref{ex:notaff} that $N$-nicely affine designs may have blocks of different sizes.

\begin{definition}{\rm
A design ${\cal D}$ with subgroup $G$ of automorphisms is called {\it $G$-pairwise transitive} if  $G$ is transitive on the following six  (possibly empty) sets: 
incident and non-incident point-block pairs,  collinear and non-collinear point pairs, intersecting and non-intersecting  blocks pairs.
}\end{definition}
For example, the points and hyperplanes of a projective geometry or of an affine geometry, with inclusion for incidence, 
are  pairwise transitive. In the latter case, taking $N$ to be the group of translations, the design is also $N$-nicely affine. 
We see in Section~\ref{designs} that the graph theoretic property in Definition~\ref{def:starlike} below is equivalent to the 
existence of a normal quotient isomorphic to an \emph{$r$-star} (a complete bipartite graph $\bfK_{1,r}$ with biparts of sizes 
1 and $r$). It is defined in terms of a subgroup of the automorphism group $\Aut(\Gamma)$ of a graph $\Gamma$.

\begin{definition}\label{def:starlike}{\rm Let $\Gamma$ be a connected bipartite graph with ordered bipartition $(B|B')$. 
We say that $\Gamma$ is {\it  $r$-starlike relative to $N$} if $N\leq \Aut(\Gamma)$, $r$ is an integer, $r\geq 2$, $N$ is 
transitive on $B$, and has $r$ orbits on $B'$.
}\end{definition}

The first result of this paper is the following characterisation.

\begin{theorem}\label{main}
Let $\Gamma$ be a connected bipartite graph with ordered bipartition $(B|B')$, $G\leq\Aut(\Gamma)$, $1\ne N\lhd G$, and $r$  an integer with $r\geq 3$.
Then the following are equivalent:
\begin{enumerate}[(a)]
 \item $\Gamma$ is  $r$-starlike relative to $N$, and $\Gamma$ is locally $(G,4)$-distance transitive;
\item the adjacency design ${\cal D}(\Gamma)$ is   $G$-pairwise transitive and $N$-nicely affine with $r$ parallel classes of blocks.
\end{enumerate}
\end{theorem}

We prove in Proposition \ref{prop:affine} that for $r\geq 3$, $s$ is at most $4$ for a locally $(G,s)$-distance transitive 
$r$-starlike graph relative to a normal subgroup of $G$, so Theorem \ref{main} addresses the highest possible value for $s$.  
This is one point of difference between locally $s$-distance transitive graphs and locally $s$-arc transitive graphs: locally $(G,s)$-arc transitive graphs that are
 $r$-starlike relative to $N\lhd G$ with $r\geq 3$ have $s\leq 3$ \cite[Lemma 5.6]{GLP} (see Section \ref{graphs} for definitions.)
 Every locally $s$-arc transitive graph is locally $s$-distance transitive but the converse does not hold. Connections between starlike locally $s$-arc transitive graphs, partial linear spaces  and homogeneous factorisations were explored in \cite{GLPpls}, while the basic starlike locally $(G,s)$-arc transitive graphs outlined in the program initiated in \cite{GLP} were investigated in \cite{GLPstar}.

For a graph $\Gamma$, subgroup $G\leq\Aut(\Gamma)$,  vertex $v$, and  $G$-invariant subset $X$ of vertices or edges  of $\Gamma$, we
denote by $\Gamma(v)$ the set of vertices adjacent to $v$ in $\Gamma$, and by $G^X$ the permutation group induced by $G$ on $X$. 
If $G^X\cong G$, then the action of $G$ on $X$ is said to be \emph{faithful}. The {\it rank} of a transitive permutation group 
is the number of orbits of a point stabiliser, and groups of rank 2 are usually called \emph{$2$-transitive}. The group $G$ 
of Theorem~\ref{main} induces several transitive actions with ranks at most 3 (see Theorem~\ref{properties}). 
This suggests several directions for further research which we discuss at the end of this section. 

\begin{theorem}\label{properties}
 If $\Gamma$, $G$, $N$ satisfy  the equivalent conditions of Theorem $\ref{main}$, then $G^B\cong G$ has rank $2$ or $3$, $G^{B'}\cong G$ is imprimitive of 
rank $3$, and for $x\in B$, $G_x^{\Gamma(x)}$ is $2$-transitive of degree $r$. 
\end{theorem}

We present in Examples \ref{ex:affspace} and \ref{selfdualex} several naturally occurring families of graphs/designs satisfying the equivalent properties of Theorem \ref{main} for which the group $N$ is elementary abelian and regular on the point set $\cal P$, so that $\cal P$ has the structure of a finite vector space. Construction~\ref{constr-regN} below properly includes these examples. The point set is a finite vector space $V$, and for a subspace $M$, $V/M$ denotes the quotient space, and $M^*=M\setminus\{0\}$, $(V/M)^*=\{v+M|v\in V\setminus M\}$ denote the sets of non-trivial elements of $M, V/M$ respectively.

\begin{construction}\label{constr-regN}
Let $V=V(d,p)$ be a vector space with group of translations $N$ and let $G=N.G_0\leq \AGL(d,p)$, where $p$ is a prime and $G_0\leq \GL(d,p)$, such that the following conditions hold.
 \begin{itemize}
 \item[(a)] $G^V$ has rank $2$ or $3$;
\item[(b)] there exists a $G_0$-orbit $\calM=\{M_1,\ldots,M_r\}$ ($r\geq 3$) of subspaces of $V$ such that $G_0^{\calM}$ is $2$-transitive (not necessarily faithful);
\item[(c)] $V=M_1+M_2$;
\item[(d)] the stabiliser $(G_0)_{M_1}$ acts transitively on $(V/M_1)^*$;
\item[(e)] $\cup_{i=1}^rM_i^*$ is a $G_0$-orbit.
\end{itemize} 
Define the design ${\cal D}=(V, \cup_{i=1}^r V/M_i, {\cal I})$ with  incidence $\cal I$ given by inclusion and let $\Gamma=\Gamma(\cal D)$ be its incidence graph.
\end{construction}

\begin{theorem}\label{constr->main}
The design ${\cal D}$ and graph $\Gamma$ in Construction $\ref{constr-regN}$ admit $G$ as an automorphism group, and the equivalent conditions of Theorem $\ref{main}$ hold.
\end{theorem}

In Construction~\ref{constr-regN}, condition (c) implies that $\dim(M_1)\geq \dim(V)/2$ and condition (e) implies that 
$G^V$ has rank 2 if and only if $\cup_{i=1}^rM_i=V$. This rather general construction motivated us to look closely at 
the graphs and designs in Theorem~\ref{main} with $N$ regular on the bipart $B$, and it leads to the following characterisation.

\begin{theorem}\label{thm:regN}
Let $\Gamma, G, N, (B|B')$ satisfy Theorem $\ref{main}(a)$ and $(b)$. If $N$ is regular on $B$, then $N$ is elementary 
abelian and $\Gamma$, ${\cal D}(\Gamma)$  can be obtained from Construction $\ref{constr-regN}$.
\end{theorem}

To complete this discussion of $r$-starlike, locally $(G,s)$-distance transitive graphs $\Gamma$, we consider briefly the case $r=2$. Although $s\leq 4$ when $r\geq3$, there is no bound on $s$ when $r=2$, see Remark~\ref{r=2list}(b), but as long as some vertex has valency at least 3, the value of $s$ is at most 14, as we will show in Theorem \ref{r=2}. Moreover, if $\Gamma$ is not a complete bipartite graph, then the adjacency design of $\Gamma$ is  resolvable (as defined in Section~\ref{designs})  but it turns out to be more useful to consider $\Gamma$ as a subdivision graph $S(\Sigma)$ of a bipartite graph $\Sigma$. The {\it subdivision graph} of $\Sigma$ is the bipartite graph with ordered bipartition $(E\Sigma|V\Sigma)$ and adjacency given by containment.

\begin{theorem}\label{r=2}
Let $\Gamma$ be a  connected bipartite graph, let $G\leq \Aut(\Gamma)$, and let $s$ be an integer such that $2\leq s\leq \diam(\Gamma)$. 
If $\Gamma$ is locally $(G,s)$-distance transitive and  $2$-starlike relative to a normal subgroup $N$ of $G$, then either (i) $\Gamma$ is the  complete bipartite graph $\bfK_{n,2}$ and $s=2$, or 
(ii) there is a bipartite graph $\Sigma$ such that $\Gamma=S(\Sigma)$ and either $\Sigma$ is known explicitly 
(Remark $\ref{r=2list} (c)$)  or  $s< \diam(\Gamma)$ and $\Sigma$ is  $(G, \lceil \frac{s+1}{2}\rceil)$-arc transitive.
In all cases either $s\leq 14$ or $\Gamma=C_{4\ell}$ for some $\ell\geq s/2$. \end{theorem}
 
The following remark
considers the extent to which the graphs $K_{n,m}$, $C_{4\ell}$, and the graphs $S(\Sigma)$ have the properties of being locally $(G,s)$-distance transitive and 2-starlike. 

\begin{remark}\label{r=2list}{\rm  
(a) The complete bipartite graph $\bfK_{n,2}$ is locally $(G,2)$-distance transitive for the group $G=S_n\times S_2$ and is $2$-starlike for the normal subgroup $N=S_n\times 1$, for instance.

(b) Suppose that $\Gamma=S(\Sigma)$ for a bipartite  $G$-arc transitive graph $\Sigma$. If $\Sigma$ has valency 2 then $\Gamma$ 
is a cycle which we discuss in part (d), so suppose that $\Sigma$ has valency at least 3.  Then $\Gamma$ is $2$-starlike relative to the 
setwise stabiliser in $G$ of the $\Sigma$-biparts. The possibilities for $\Sigma$ such that $\Gamma$ is locally 
$(G,\diam(\Gamma))$-distance transitive, for some $G$, can be determined from the classification in \cite{DD}. They are the graphs 
 $\bfK_{n,n}$ ($n\geq 3$),  and the incidence graphs of the following generalised polygons: Desarguesian
projective planes,   symplectic generalised quadrangles over 
finite fields of characteristic 2, and  split-Cayley generalised hexagons over  finite fields of characteristic 3.
Details about the groups $G$ can be found in \cite{DD}. In particular, for every example, there exists an automorphism group $G$ such that $\Gamma$ is locally $(G,\diam(\Gamma))$-distance transitive and $\Gamma$ is 2-starlike relative to a normal subgroup of $G$ (see Lemma \ref{lem:last}).

(c) If $\Gamma=S(\Sigma)$ and $\Sigma$ is  $(G, \lceil \frac{s+1}{2}\rceil)$-arc transitive for some   $s< \diam(\Gamma)$, then $\Gamma$ is locally $(G,\diam(\Gamma))$-distance transitive and $\Gamma$ is 2-starlike relative to a normal subgroup of $G$ (see Lemma \ref{lem:last}).

(d) The cycle $\Gamma=C_{4\ell}$ is the subdivision graph of the smaller bipartite cycle $\Sigma = C_{2\ell}$.  The group $G=\Aut(\Sigma)\cong D_{4\ell}$ has a normal subgroup $N\cong D_{2\ell}$ which is transitive on $E\Sigma$ and has two orbits on $V\Sigma$, so that $\Gamma$ is locally $(G,s)$-distance transitive for all $s\leq 2\ell$, and is $2$-starlike relative to $N$. 

(e) Each complete bipartite graph $\Gamma=\bfK_{n,m}$  has diameter 2, is  locally $(G,2)$-distance 
transitive (for example with $G=S_n\times S_m$),  and is $2$-starlike (for example with respect to $N=S_n\times 
(S_{k}\times S_{m-k})$). However,  we see in the proof of Theorem \ref{r=2} that for $m\geq 3$,  there is no locally 2-distance transitive group $G$ with normal subgroup $N$ such that $\Gamma$ is 2-starlike relative to $N$. However, for $m$ even,  there is a group $G$ with normal subgroup $N$ such that $\Gamma$ is locally $(G,1)$-distance transitive and 2-starlike relative to $N$: take for instance $G=S_n\times (S_{m/2}\wr S_2)$  and $N=S_n\times (S_{m/2}\times S_{m/2})$.
} \end{remark}

In Section 2, we explain some design theoretic and graph theoretic concepts and prove some preliminary results on graphs, designs and the links between the two (for affine designs).
In Section 3 we display some examples and prove Theorem \ref{constr->main}. In Section 4, we prove Theorems \ref{main} and \ref{properties}, and compute the intersection arrays of a graph satisfying the conditions of Theorem \ref{main}. In Section 5, we study the case where $N$ has a regular action on points, and finally in Section 6 we prove Theorem \ref{r=2}.

\subsection*{Commentary and future directions}\label{openprobs} 
The information in Theorem~\ref{properties} suggests directions for 
further  research (some of which we intend to pursue) since the finite 2-transitive permutation groups, and quasiprimitive 
rank 3 groups are known explicitly as a consequence of the finite simple group classification \cite{Cam99,DGLPP,Lieb,LiebSa}. 
The most general (and probably very difficult) problem is the following.

\begin{problem}\label{pwtr}{\rm Classify finite pairwise transitive designs. 
}\end{problem}

A particularly interesting subfamily are the pairwise transitive $2$-designs where each point pair lies in at 
least one common block. These designs are symmetric if each block pair intersects nontrivially, and otherwise 
they are quasisymmetric (where block pairs have two possible intersection sizes). Both types of 2-designs have 
been studied extensively, but without the assumption of pairwise transitivity, see \cite{DesignTheory,Lander,quasisym}. We plan to exploit the information in Theorems~\ref{properties} and~\ref{thm:regN} to study pairwise transitive 2-designs. 

However the general case of Problem~\ref{pwtr} remains completely open. A second special case which is important for the application to starlike locally $s$-distance transitive graphs is the following.

\begin{problem}\label{naff}{\rm Classify $G$-pairwise transitive, $N$-nicely affine designs, where $N\lhd G$. 
}\end{problem}

Note that a complete solution to Problem~\ref{naff} would give a classification of locally $(G,4)$-distance transitive $r$-starlike graphs. This problem is beyond our reach at present. We intend to study the special case in which the group $G$ is quasiprimitive on points (each nontrivial normal subgroup transitive). A solution of this special case will yield a classification of  the $G$-basic locally $(G,4)$-distance transitive $r$-starlike graphs (identified as important in \cite{DGLP}, see also Section~\ref{graphs}).

\section{Definitions, examples and preliminary results}\label{Def}

\subsection{Design theoretic concepts}\label{designs}  

A design ${\cal D}=({\cal P}, {\cal B}, {\cal I})$ consists of a point set ${\cal P}$, a block set ${\cal B}$ 
and an incidence relation ${\cal I}\subseteq {\cal P} \times {\cal B}$. The relation 
${\cal I}$ induces the following relations on ordered pairs of distinct objects of ${\cal D}$: a point-block 
pair is either {\it incident} if it lies in ${\cal I}$, or {\it non-incident} if it does not; a point pair is 
{\it collinear} if the two points are incident with at least one common block, and otherwise is {\it non-collinear}; 
a block pair is {\it intersecting} if the two blocks are incident with at least one common point, and {\it non-intersecting} 
otherwise.
An \emph{automorphism} of a design is a permutation of ${\cal P}\cup {\cal B}$ preserving points, blocks, and incidence.
We say that  ${\cal D}$ has \emph{no repeated blocks} if there are no two distinct blocks incident with exactly the same point-sets. By proposition \ref{star quotient}, we see that the designs with the properties of Theorem \ref{main} have no repeated blocks.

 The \emph{dual design}  is the design ${\cal D}^*=({\cal B}, {\cal P}, {\cal I}^*)$, where ${\cal I}^*
=\{(b,x)\,|\,(x,b)\in{\cal I}\}$. A design $\cal D$ is {\it connected} if its incidence graph $\Gamma({\cal D})$ is connected.
For connected designs the set of intersecting block pairs is always non-empty.

A design $\cal D$ is called a {\it $t-(v,k,\lambda)$-design} (see \cite{DesignTheory} for instance) if $|{\cal P}|=v$, 
each block is incident with $k$ points, and each $t$-subset of points is incident with exactly $\lambda$ blocks; $\cal D$ 
is {\it non-trivial} if $1\leq t<k<v$. 
A {\it $t$-design} is a  $t-(v,k,\lambda)$-design for some parameters $v,k,\lambda$. 

A design whose automorphism group is transitive on points and blocks is automatically a $1$-design, and so is its dual design. 
In particular, a $G$-pairwise transitive design $\cal D$ is point-transitive and block-transitive, and so is a $1$-design. 
Moreover if each point pair is collinear then $G$ is $2$-transitive on points, and so $\cal D$ is a $2$-design. 
If $\cal D$ is connected and pairwise transitive, and its set of non-intersecting block pairs  is non-empty,  
then there are exactly two possible intersection sizes for block pairs (one of them being $0$), and 
the design is {\it quasisymmetric}. In particular, pairwise transitive resolvable designs are quasisymmetric (see below).

Recall the concept of a parallel class, explained after Definition \ref{def:nice}. 
A design ${\cal D}$ is {\it resolvable} if its block set admits a partition into parallel classes.
An {\it affine $t$-design} (see \cite{DesignTheory} for instance) is a resolvable $t$-design for which there is a positive 
constant $\mu$ such that any two blocks in distinct parallel classes are incident with exactly $\mu$ common points.  
In particular an $N$-nicely affine design with blocks of a constant size is an affine 1-design. As we see in 
Example~\ref{ex:degen}, there are some degenerate 
disconnected examples which have a single parallel class of blocks.  We also give, in Example~\ref{ex:notaff} 
the promised examples of  $N$-nicely affine designs with blocks of different sizes.

\begin{example}\label{ex:degen}{\rm 
Let $k,\ell$ be positive integers with $\ell\geq2$, let $X$ be a set of size $k\ell$, and let $\cal B$ be a partition of 
$X$ with $\ell$ parts of size $k$. Let $G=N=S_k\wr S_\ell$ denote the stabiliser in $S_{k\ell}$ of the partition $\cal B$,
and let ${\cal D} = (X,{\cal B}, {\cal I})$ with $\cal I$  natural inclusion. Then $\cal D$ is $G$-pairwise transitive (notice the set of pairs of intersecting lines is empty, so there is only 5 transitivity properties to check) and $N$-nicely affine (with only one  parallel class of blocks). Its incidence graph $\Gamma({\cal D})= \ell.\bfK_{1,k}$.}
\end{example}
Moreover we have the following lemma.
\begin{lemma}
\begin{enumerate}
 \item[(i)] Let $\cal D$ be a disconnected pairwise transitive design with no repeated blocks. Then $\cal D$ is as in  Example $\ref{ex:degen}$. 
\item[(ii)] Let $\cal D$ be a  nicely affine designs with a unique parallel class of blocks and no repeated blocks. Then $\cal D$ is as in  Example $\ref{ex:degen}$. 
\end{enumerate}
\end{lemma}
\proof
(i) Let $\cal D$ be a disconnected pairwise transitive design with no repeated blocks. Let $k$ be the number of points in each component and let $\ell$ be the number of components. Since  $\cal D$  is pairwise transitive, its automorphism group is in particular point-transitive, and so all the connected components are isomorphic. Let $x$ be a point in one component and $L$ be a line in another component. Obviously $x$ and $L$ are non-incident. Since the automorphsim group is transitive on the set of non-incident point-block pairs, any two non-incident point and block are in distinct components. So all the blocks in one component are incident with all the points in that component. Since $\cal D$ has no repeated blocks, it follows that  $\cal D$ is as in  Example \ref{ex:degen}.\\
(ii)  Let $\cal D$ be a  nicely affine designs with a unique parallel class of blocks. Then the blocks of $\cal D$ are pairwise non-intersecting, and so  $\cal D$ is a disjoint union of lines. Since there is a group $N$ transitive on points, all the blocks must be incident with the same number of points, and $\cal D$ is as in  Example \ref{ex:degen}.
\qed

\begin{example}\label{ex:notaff}{\rm 
Let $k,\ell$ be positive integers such that $k>\ell>1$, and let $N=K\times L\cong Z_k\times Z_\ell$. Define the 
design $\cal D$ with point set $N$, block set consisting of the subsets $b_i:=\{(i,j)\,|\,j\in L\}$ for $i\in K$, 
and $c_j:=\{(i,j)\,|\,i\in K\}$ for $j\in L$, and inclusion as incidence. Then $|b_i|=\ell$ and $|c_j|=k$, $N$ acts transitively on points by multiplication, and the two $N$-orbits on blocks are parallel classes (namely the $b_i$ and the $c_j$). Since $b_i\cap c_j=\{(i,j)\}$, $\cal D$ is $N$-nicely affine.
}
\end{example}

Pairwise transitivity gives strong restrictions on the actions of the automorphism group on points and blocks. 

\begin{lemma}\label{lem:pt}
 Let  ${\cal D}$ be a $G$-pairwise transitive design. 
Then $G$ has rank $2$ or $3$ on points and rank $2$ or $3$ on blocks.
Moreover the following statements hold.
\begin{itemize}
 \item[(a)] If  ${\cal D}$ is resolvable then $G$ has rank $3$ and is imprimitive on blocks.
\item[(b)] If  ${\cal D}$ is $N$-nicely affine for some $N\leq G$, then for any point $x$, $G_x$ is $2$-transitive on the set of $N$-orbits on blocks.
\end{itemize}
\end{lemma}
\proof 
A point stabiliser in $G$ has at most 3 orbits on points: itself, the points collinear to it, and the points non-collinear to it (this last set could be empty).
Similarly   a block stabiliser has at most 3 orbits on blocks: itself, the blocks intersecting it, and the blocks not intersecting it (this last set could be empty). 
Hence the first statement holds.

(a) Assume that ${\cal D}$ is resolvable. Then the set of blocks not intersecting a given block is non-empty (it contains the parallel blocks), so $G$ has rank 3 on blocks. Moreover the parallel classes are preserved by $G$, and so they form a system of imprimitivity.

(b)  Assume that ${\cal D}$ is $N$-nicely affine (so in particular it is resolvable). Since two blocks in different parallel classes do intersect, $G$ is 2-transitive in its action on the parallel classes. Suppose $g\in G$ maps two given parallel classes to two given parallel classes. Since $N$ is transitive on points, there is an element $n\in N$    mapping $x^g$ to $x$, so that $gn$ is in $G_x$. By definition, the set of $N$-orbits are the parallel classes, so $gn\in G$ has the same action on the parallel classes as $g$. It follows that  $G_x$ is 2-transitive on the parallel classes.
\qed

\subsection{Graph theoretic concepts}\label{graphs}  

A graph $\Gamma$ consists of a vertex set $V\Gamma$ and a subset $E\Gamma$of unordered pairs from $V\Gamma$, called edges. If 
$e=\{v,u\}\in E\Gamma$, we say that $v$ is incident with $e$ and adjacent with $u$. 
For positive integers $i\leq \diam(\Gamma)$, we denote by $\Gamma_i(x)$ the set of vertices at 
distance $i$ from $x$ in $\Gamma$ (so $\Gamma(x)=\Gamma_1(x)$). 
Recall that a locally $(G,s)$-distance transitive graph is a graph with automorphism group $G$ such that for every vertex $x$, $G_x$ is transitive on $\Gamma_i(x)$ for each $i\leq s$.

An \emph{$s$-arc} starting at $v_0$ in $\Gamma$ is an $(s+1)$-tuple $v_0,v_1,\ldots,v_s$ of vertices such that consecutive vertices are adjacent and $v_{i-1}\neq v_{i+1}$ for $i=1,2,\ldots,s-1$. 
A \emph{locally $(G,s)$-arc transitive graph} is a graph with automorphism group $G$ such that for every vertex $x$, $G_x$ is transitive on the set of $s$-arcs starting at $v$. If moreover $G$ is transitive on vertices then the graph is said to be \emph{ $(G,s)$-arc transitive}.

For a bipartite graph $\Gamma$ with 
ordered bipartition $(B|B')$, each edge meets each of $B$ and $B'$ in a single vertex, and we sometimes identify 
$E\Gamma$ with the corresponding subset $\cal I$ of $B\times B'$, (note $\cal I$ is the incidence relation of the 
adjacency design ${\cal D}(\Gamma)$). We denote by $\bfK_{m,n}$ the complete bipartite graph with biparts of sizes $m$ and $n$.

If $N\lhd G\leq \Aut(\Gamma)$, the {\it $G$-normal quotient} of $\Gamma$ relative to $N$ is the graph $\Gamma_N$  whose vertices are the  $N$-orbits in $V\Gamma$, such that two $N$-orbits  are adjacent if and only if there exist two vertices, one in each of the $N$-orbits, forming an edge of $\Gamma$. 
We say that a $G$-edge-transitive graph $\Gamma$ is \emph{$G$-basic} if each $G$-normal quotient of $\Gamma$ is one of $\bfK_1$, $\bfK_2$ or an $r$-star $\bfK_{1,r}$ for some $r\geq2$. Edge-transitive graphs with star normal quotients are starlike, as we now show.

\begin{proposition}\label{starquo=starlike}
Let $\Gamma$ be a connected $G$-edge-transitive graph and let $N\lhd G\leq \Aut(\Gamma)$.
Then $\Gamma_N$ is isomorphic to $\bfK_{1,r}$ if and only if $\Gamma$ is $r$-starlike relative to $N$. 
\end{proposition}
\proof
Assume $\Gamma_N\cong\bfK_{1,r}$. Then $N$ has vertex orbits $v^N$, $x_i^N$ ($1\leq i\leq r$), there is at least one edge of $\Gamma$ between  the sets $v^N$ and $x_i^N$ for each $i$,  there is no edge between  $x_i^N$ and $x_j^N$ (for $i\neq j$), and since $G$ is transitive on $E\Gamma$, there is no edge between vertices in the same $N$-orbit. Thus $\Gamma$ is bipartite with ordered bipartition $(v^N| \cup_{i=1}^r x_i^N)$ and is $r$-starlike relative to $N$.  

Conversely, assume  $\Gamma$ is $r$-starlike relative to $N$, with bipartition $(B|B')$. Then there are no edges among the vertices of $B$, nor among the vertices of $B'$. Since $\Gamma$ is connected and $G$-edge-transitive, there must be an edge between $B$ and each $N$-orbit in $B'$. Thus $\Gamma_N\cong\bfK_{1,r}$. 
\qed

\begin{lemma}\label{lem:diffneigh}
 Assume $\Gamma$ is $r$-starlike ($r\geq2$) relative to $N$ with ordered bipartition $(B|B')$, and $\Gamma$ is  locally $(G,2)$-distance transitive with  $1\neq N\lhd G$.
Then either (i) $\Gamma$ is complete bipartite, or (ii) $\Gamma(x)=\Gamma(x')$ implies $x=x'$, and $G$ acts faithfully on both $B$ and $B'$.
\end{lemma}
\proof
Suppose that  $\Gamma(x)=\Gamma(x')$ and $x\ne x'$. Since $\Gamma$ is bipartite, $x$ and $x'$ must be in the same bipart, say $B$. Then $\dist_\Gamma(x,x')=2$. Since $G_x$ is transitive on $\Gamma_2(x)$, we have 
$\Gamma(x)=\Gamma(y)$ for all $y\in\Gamma_2(x)$. Repeating this argument, since $\Gamma$ is connected, we find that $\Gamma(x)=\Gamma(y)$ 
for all $y\in B$. Thus  $\Gamma$ is complete bipartite. 
\qed

\subsection{Links between designs and graphs: the affine case}

For affine designs, there is a direct correspondence between pairwise transitivity of $\calD(\Gamma)$ and local 4-distance transitivity 
of $\Gamma$. Note that, if ${\cal D}(\Gamma)$ is affine, then $\Gamma$ has diameter at most $4$.

\begin{proposition}\label{4dist=PT}
Suppose that $\Gamma$ is a connected bipartite graph with ordered bipartition $(B|B')$  
such that ${\cal D}(\Gamma)$ is affine. 
Let $G$ be an automorphism group of $\Gamma$ (and hence also of $\mathcal{D}(\Gamma)$)
fixing $B, B'$ setwise, and let $x\in V\Gamma$.
\begin{enumerate}
   \item[(i)] If $G$ is transitive on $B$ and $B'$ then, for each line of Table~$\ref{tbl:equiv}$, $G$ is transitive on the set of 
pairs of elements of $\calD(\Gamma)$ described in column $1$ if and only if $x$ lies in the bipart in column $2$ and 
$G_x$ is transitive on the set $\Gamma_i(x)$  described in  column $3$.
\item[(ii)] The graph $\Gamma$ is locally $(G,4)$-distance transitive if and only if $\mathcal{D}(\Gamma)$ is $G$-pairwise transitive.
\end{enumerate}
\end{proposition}

\begin{table}[ht]
\begin{center}
\begin{tabular}{|l|c|l|}
\hline 
Pairs in $\mathcal{D}(\Gamma)$ & $x$ in &Subset of $V\Gamma$ \\
\hline 
incident point-block pairs & $B$ & $\Gamma_1(x)$  \\
incident point-block pairs & $B'$&$\Gamma_1(x)$ \\
collinear point pairs &$B$& $\Gamma_2(x)$ \\
intersecting block pairs &$B'$& $\Gamma_2(x)$   \\
non-incident point-block pairs &$B$& $\Gamma_3(x)$  \\
non-incident point-block pairs &$B'$& $\Gamma_3(x)$  \\
non-collinear point pairs &$B$& $\Gamma_4(x)$  \\
non-intersecting block pairs &$B'$& $\Gamma_4(x)$ \\
\hline
\end{tabular}
\caption{Design and graph symmetry properties}\label{tbl:equiv}
\end{center}
\end{table}

\proof
The proof of part (i) follows immediately from the definitions of $\Gamma_i(x)$ and $\calD(\Gamma)$,
and we note that for the last four lines we use the fact that $\cal D(\Gamma)$ is affine.
Also,  the design is $G$-pairwise transitive if and only if 
$G$ is transitive on the sets of pairs in each of the lines of Table~\ref{tbl:equiv}, and 
if $G$ is transitive on $B$ and $B'$, then the graph is locally 
$(G,4)$-distance transitive if and only if $G_x$ is transitive on the set $\Gamma_i(x)$ for each of the lines of
the table.
Notice that each of the conditions ``$\Gamma$ is locally $(G,4)$-distance transitive'' and ``$\mathcal{D}(\Gamma)$  
is $G$-pairwise transitive'' implies that $G$ is transitive on $B$ (points) and on $B'$ (lines).
Thus part (ii) follows from part (i).
\qed

\section{Examples and Proof of Theorem \ref{constr->main}}

We present three families of examples of pairwise transitive designs, some of which are nicely affine. It will follow from Theorem \ref{main} that the incidence graphs of the nicely affine examples are $r$-starlike and locally $4$-distance transitive. Further, the incidence graphs of the examples in the first two families are basic, but those in the third family are not.

\begin{example}\label{ex:pt1}
{\rm Let $X$ be a set of $v$ points ($v\geq3$), let $X^{\{v-1\}}$ denote the set of $(v-1)$-subsets of $X$, and let $\cal I$ be the inclusion relation. Then ${\cal D}:=(X, X^{\{v-1\}}, {\cal I})$ is  a $1-(v,v-1,v-2)$ design and also a $2-(v,v-1,v-3)$ design, with automorphism group $G=S_v$, the symmetric group on $X$. This design is $G$-pairwise transitive, and its incidence graph $\Gamma(\mathcal{D})$ is $G$-basic since its only proper $G$-normal quotient is $\bfK_2$ (relative to $N=G$ or $N=A_v$, both of which are transitive on $X$ and $X^{\{v-1\}}$).
}
\end{example}

Proving pairwise transitivity in Example~\ref{ex:pt1} requires an easy check of $G$-transitivity on four sets of pairs: incident point-block pairs, non-incident point-block pairs,  collinear point pairs, and intersecting block pairs. Note that there are no 
 non-intersecting block pairs and no non-collinear point pairs in these examples.

\begin{example}\label{ex:affspace}{\rm
Let $\cal P, \cal B$ be the sets of  points and hyperplanes of an affine space $AG(d,q)$ (where $d\geq3$), and let $\cal I$ be the inclusion relation. Then the design ${\cal D}:=({\cal P}, {\cal B}, {\cal I})$ is  a $1-(q^d,q^{d-1},\frac{q^d-1}{q-1})$ design and also a $2-(q^d,q^{d-1},\frac{q^{d-1}-1}{q-1})$ design. Let $G=AGL(d,q)$, and let $N$ be the group of translations. Then $\cal D$ is  
$N$-nicely affine with intersection size $\mu=q^{d-2}$, and $G$-pairwise transitive, and 
the incidence graph $\Gamma(\mathcal{D})$ is $G$-basic. Moreover, the proper $G$-normal quotients are the star $\bfK_{1,r}$ with $r=\frac{q^d-1}{q-1}$, and $\bfK_2$ (see below).
}
\end{example}

The design $\cal D$ in Example~\ref{ex:affspace} is $N$-nicely affine since $N$ preserves each parallel class of hyperplanes and any two non-parallel hyperplanes intersect in $\mu=q^{d-2}$ points. Proving pairwise transitivity requires checking transitivity on the five sets: incident and non-incident point-block pairs,  collinear point pairs, and intersecting and non-intersecting block pairs. This is easy using the fact that $G_0=\GL(d,q)$ is transitive on the ordered bases of $V(d,q)$.
Finally since $N$, the unique minimal normal subgroup of $G$, is point-transitive and transitive on each parallel class of blocks, and since each normal subgroup properly containing $G$ is block transitive, it follows that the proper $G$-normal quotients are the star $\bfK_{1,r}$ with $r=\frac{q^d-1}{q-1}$, and $\bfK_2$. Thus $\cal D$ is $G$-basic.

\begin{example}\label{selfdualex}{\rm
Let $\cal D$, $G, N$ be as in Example~\ref{ex:affspace}, let $u\in{\cal P}\setminus\{0\}$, and let ${\cal B}'$ the 
subset of $\cal B$ consisting of hyperplanes not containing a line with direction $\langle u\rangle$.  Let ${\cal I}'$ be the restriction of $\cal I$ 
to ${\cal P}\times {\cal B}'$, and let $K=NU$, where  $U$ is the stabiliser in $\GL(d,q)$ of $\langle u\rangle$. Then 
${\cal D}'=({\cal P}, {\cal B}',{\cal I}')$ is a self-dual $1-(q^d,q^{d-1},q^{d-1})$ design, which is $K$-pairwise 
transitive and $N$-nicely affine, but is not $K$-basic.   
}
\end{example}

\begin{lemma}
The assertions of Example $\ref{selfdualex}$ are valid.
\end{lemma}
\proof
The design $\calD'$ is $N$-nicely affine since the $N$-orbits on blocks are the parallel classes of hyperplanes in $\cal B'$, $N$ is transitive 
on points, and any two non-parallel blocks intersect in $q^{d-2}$ points.
Let $M\leq N$ consist of the translations stabilising the 1-dimensional subspace $\langle u\rangle$. Then $M\lhd K$. 
The $M$-orbits in $\cal P$ are the lines parallel to $\langle u\rangle$, and the $M$-orbits in ${\cal B}'$ are 
parallel classes of ${\cal B}'$, so the quotient $\Gamma(\mathcal{D})_M$ is none of $\bfK_1$, $\bfK_2$ or a star. Thus
$\mathcal{D}'$ is not $K$-basic. 

Next we display an automorphism of $\Gamma(\mathcal{D})$ switching $\cal P$ and ${\cal B}'$. 
Replacing $\calD'$ by an isomorphic image under some element of $G$, we may assume that $u=(1,0,\dots,0)$.
For a subset $S\subseteq\calP$, we denote by $S^\perp$ the set of points $v$ such that $v\cdot x=0$ for all $x\in S$
(where `$\cdot$' denotes the usual inner product). In particular 
we let $U_0= \langle u\rangle^{\perp}$.
Then each point is of the form $au+v_0$ for some $v_0\in U_0$ amd $a\in F_q$, and each hyperplane in ${\cal B}'$ 
is of the form $w+\langle d\rangle^{\perp}$ where $w\in\calP$, and the normal vector $d\not\in U_0$ so $d$ can be chosen 
such that $d\cdot u=1$.
Define $\phi:{\cal P}\rightarrow {\cal B}'$ by $au+v_0\mapsto -au+\langle u+v_0\rangle^{\perp}$, and $\theta: {\cal B}'
\rightarrow \cal P$ by $w+\langle d\rangle^{\perp}\mapsto -(w\cdot d)u+(d-u)$. 
The map $\theta$ is well-defined since, if $w_1+\langle d\rangle^{\perp}=w_2+\langle d\rangle^{\perp}$, then 
$w_1-w_2\in \langle d\rangle^{\perp}$ and so $w_1\cdot d=w_2\cdot d$.
Observe that $\phi(au+v_0)\in {\cal B}'$: indeed $u+v_0\notin U_0$ since $(u+v_0).\cdot u=u\cdot u=1$.
It is an easy exercise to check that $\phi({\cal P})= {\cal B}'$ and $\theta=\phi^{-1}$. 
We claim that the permutation $\alpha$ of the vertex set ${\cal P}\cup{\cal B}'$ of $\Gamma(\mathcal{D})$, 
defined by $\alpha|_{\cal P}=\phi$ and  $\alpha|_{\cal B'}=\theta$, is an automorphism of $\Gamma(\mathcal{D})$.
It is sufficient to prove that, for a hyperplane $H\in  {\cal B}'$, $v\in H$ implies that $H^{\theta}\in v^{\phi}$.
Suppose that $H=v+\langle d\rangle^{\perp}$ (with $u\cdot d=1$) and $v=au+v_0\in H$. 
Then $H^{\theta}=-(v\cdot d)u+(d-u)=-((au+v_0)\cdot d)u+(d-u)=-au-(v_0\cdot d+1)u+d$ and $v^{\phi}=-au+\langle u+v_0\rangle^{\perp}$.
Therefore $H^{\theta}\in v^{\phi}$ if and only if $H^{\theta}+au\in \langle u+v_0\rangle^{\perp}$.
Now $(H^{\theta}+au)\cdot(u+v_0)=(-(v_0\cdot d+1)u+d)\cdot(u+v_0)=-(v_0\cdot d+1)+1+v_0\cdot d=0$, and 
hence $H^{\theta}\in v^{\phi}$. This proves the claim, and hence we have shown that $\mathcal{D}'$ is self-dual.

The fact that  $\mathcal{D}'$ is self-dual simplifies the proof that it is $K$-pairwise transitive. 
We only need to check transitivity on  4 sets: incident and non-incident point-block pairs, and collinear and 
non-collinear point pairs. This is easy using the fact that $K_0=U$ is transitive on the ordered bases of 
$V(d,q)$ with first basis element in $\langle u\rangle$.
\qed

The last two examples can be obtained from Construction \ref{constr-regN}: for Example \ref{ex:affspace} 
we take $\calM$ to be the set of all $(d-1)$-subspaces of $V=V(d,p)$, while for  Example \ref{selfdualex} 
we take $\calM$ to be the set of all dimension $(d-1)$-subspaces of $V=V(d,p)$ not containing the vector $u$.

\medskip

\noindent
{\it Proof of Theorem $\ref{constr->main}$.} 
Let $G, V, \Gamma, \calD$ be as in Construction~\ref{constr-regN}. Then $G$ acts faithfully on $V$, and $G$ 
also has a natural induced action on the cosets of subspaces in $\calM$, and hence on the blocks of $\calD$. 
This action preserves inclusion, and hence preserves incidence in ${\cal D}$. Thus $G\leq \Aut({\cal D})$.
Now  $\calD(\Gamma)={\cal D}(\Gamma({\cal D}))={\cal D}$. For $v\in V$ let $t_v:V\rightarrow V:x\mapsto v+x$. Then 
$N=\{t_v|v\in V\}$ is the group of translations: $N$ is regular on $V$ and normal in $G$.

By conditions (c) and (b), distinct cosets $x+M_i$ and $y+M_j$ intersect non-trivially if and only if $i\neq j$,
and in this case $(x+M_i)\cap(y+M_j)=(x+y)+ (M_i\cap M_j)$, of size $|M_i\cap M_j|$, which is independent of 
$i, j$ by condition (b). Thus intersecting block pairs are coset pairs of distinct subspaces in $\cal M$, and so
${\cal D}$ is a resolvable design, with parallel classes the cosets of a fixed $M_i$. 
A translation $t_v\in N$ maps $x+M_i$ to $(x+v)+M_i$, and hence the $N$-orbits on blocks are the parallel classes. 
This implies that $\calD$ is $N$-nicely affine. It also implies, together with condition (b), that $G$ is 
transitive on blocks as well as on points. This is useful in our proof of pairwise transitivity.

By (e), $G_0$ is transitive on the points collinear with $0$ (since $\cal M$ is the set of blocks containing $0$),
and by (a),  $G_0$ is also transitive on the points non-collinear with $0$ (if any). Thus $G$ is transitive on 
collinear and non-collinear point pairs.
By (b),   $G_0$ is transitive on the blocks incident with $0$, and by (b) and (d),   $G_0$ is transitive on 
the blocks non-incident with $0$. This gives $G$-transitivity on incident and non-incident point-block pairs.

By (d), $G_{M_1}$ is transitive on blocks not intersecting (that is, parallel to)  $M_1$. Finally we
prove that $G_{M_1}$ is transitive on blocks intersecting $M_1$. This will complete the proof of pairwise transitivity.
An arbitrary block intersecting $M_1$ has the form $x+M_i$ with $i>1$. It is sufficient to map $x+M_i$ to $M_2$ by an 
element of $G_{M_1}$. By (b), $G_{M_1}$ is transitive on ${\cal M}\setminus\{M_1\}$, so there exists 
$g\in G_{M_1}$ such that $M_i^g=M_2$. By (c), $x^g=x_1+x_2$ for some 
$x_1\in M_1, x_2\in M_2$, so $(x+M_i)^g=x^g+M_2=x_1+M_2$.  Finally $G_{M_1}$ contains the subgroup  
$N_{M_1}=\{t_y|y\in M_1\}$, and $(x_1+M_2)^{t_{-x_1}}=M_2$. 
\qed

\section{Proofs of Theorems \ref{main} and \ref{properties}}

Lemma~\ref{lem:diffneigh} illustrates the exceptional role of complete bipartite graphs, 
and so we make the following hypothesis, which we assume throughout this section.

\begin{hypothesis}\label{hyp}{\rm
Let $r\geq2$, $s\geq2$. Let $\Gamma$ be an $r$-starlike graph relative to $N$ with ordered bipartition $(B|B')$ such that $\Gamma$ 
is not complete bipartite.
Suppose that $\Gamma$ is locally $(G,s)$-distance transitive and that $1\neq N\lhd G$.
Let ${\cal N}=\{B,C_1,C_2,\ldots C_r\}$ be the set of $N$-orbits on $V\Gamma$. 
}\end{hypothesis}

As in \cite{DGLP}, we define $\dist_\Gamma(S) = \min\{\dist_\Gamma(v, w)|v, w\in S, v\neq w\}$ for a subset $S$ of at 
least two vertices of $\Gamma$.  When $\Gamma$ is bipartite, we denote by $G^+$ the stabiliser of the biparts in the 
automorphism group $G$. 
The following lemma follows easily from \cite[Lemma 5.2]{DGLP}. 
\begin{lemma}\label{cor:star}
Assume that Hypothesis $\ref{hyp}$ holds.
Then $G=G^+$, and there is an integer $\ell>1$ such that, for each $C\in{\cal N}\setminus\{B\}$, $|C|=\ell$ and $\dist_\Gamma(C)\geq 4$. In particular, for each $x\in B$, $|\Gamma(x)\cap C|=1$, so $x$ has valency $r$.
\end{lemma}
\proof
Since $N\lhd G$, $G$ permutes the $N$-orbits and since $\Gamma_N$ is isomorphic to $\bfK_{1,r}$ by Proposition \ref{starquo=starlike}, we have that $G=G^+$ fixes $B$ and $B'$.
Since $\Gamma$ is $(G,1)$-distance transitive, $G$ is transitive on $E\Gamma$, and so the orbits of $G$ on vertices are precisely $B$ and $B'$. In particular, all the $N$-orbits $C_i$ have the same size, $\ell$ say. 

Since $\Gamma$ is not complete bipartite, in particular $\Gamma$ is not a star, and so  $|B|\geq 2$. 
If  $\ell=1$, then $N$ fixes $ B'$ vertexwise and is transitive on $B$, and so $\Gamma$ is complete bipartite, which we have assumed is not the case here. Thus $\ell\geq 2$.
By Lemma 5.2 of \cite{DGLP}, we have  that $\dist_\Gamma(C_i)\ge4$ for each $i$.
It follows that, for $x$ a vertex in $B$, $|\Gamma(x)\cap C_i|\leq 1$ for all $i$.
Since $\{B,C_i\}$ is an edge of $\Gamma_N$, there is a vertex in $B$ adjacent to a vertex in $C_i$.
 Since $B$ and the $C_i$'s are $N$-orbits, it follows that $|\Gamma_1(u)\cap C_i|=1$ for each $i$. Therefore the valency of the vertices in $B$ is $r$. 
\qed

We explore additional properties of star-like graphs. 

\begin{proposition}\label{star quotient}
Assume Hypothesis $\ref{hyp}$ holds.
Then  the adjacency design $\mathcal{D}(\Gamma)$ has no repeated blocks is resolvable with $r$ parallel classes.
\end{proposition}

\proof
By Lemma \ref{cor:star},  ${\cal N}=\{B,C_1,C_2,\ldots C_r\}$, where $B$ is a bipart,  $|C_i|=\ell\neq 1$ does not depend on $i$, $\dist_\Gamma(C_i)\ge4$ for all $i$, and  the valency of the vertices in $B$ is $r$. Moreover the vertices in $ B'$ all have the same valency, say $k$, and all their neighbours are in $B$.

Counting the edges between $B$ and any $C_i$ yields that $|B|=k|C_i|=k\ell$, and for each $i$, $B$ is partitioned by $\{\Gamma(w')|w'\in C_i\}$.
If $k=1$, then $\Gamma$ is the disjoint union of $\ell$ stars $\bfK_{1,r}$ and $\Gamma$ is not connected, a contradiction.  Thus $k\geq 2$ and we recall that $\ell\geq 2$.

Consider the adjacency design $\mathcal{D}(\Gamma)$. 
From Lemma \ref{lem:diffneigh}, we get that two distinct `blocks' are incident with distinct `point'-sets (hence $\mathcal{D}(\Gamma)$ has no repeated blocks) and two distinct points are incident with distinct `block'-sets. Each `block' is incident with $k$ `points' and each `point' is incident with $r$ `blocks'. So $\mathcal{D}(\Gamma)$ is a $1-(k\ell,k,r)$ design.
 Moreover, we have seen that every `point' is incident with exactly one `block' in each $C_i$. Therefore the blocks $C_1,C_2,\dots,C_r$ of ${\cal N}$ give a resolution of the `blocks' of $\mathcal{D}(\Gamma)$ into parallel classes. 
\qed

\begin{remark}{\rm
In the case  $\diam_\Gamma(v)\leq 3$ for $v\in B$, any two vertices in $B$ are at distance 2, and so, by local 
$(G,2)$-distance transitivity, the corresponding two points of $\calD(\Gamma)$ are incident to a constant number 
of blocks (in $ B'$). Therefore in that case $\mathcal{D}(\Gamma)$ is actually a 2-design. More precisely,  it is  
a $2-(k\ell,k,r(k-1)/(k\ell-1))$ design which is {quasisymmetric} with intersection numbers $0$ and $r(k-1)/(k\ell-1)$
(the possble sizes of block intersections).
The last parameter is found by counting in two different ways the number of triples $(v_1,v_2,x)$ such that 
$v_1,v_2\in B$, $x\in B'$ and $v_1,v_2$ are both adjacent to $x$.
}
\end{remark}

\begin{proposition}\label{prop:affine}
Assume Hypothesis $\ref{hyp}$ holds and that $r\geq 3$.
Then $s\le 4$. Moreover, if  $s=4$, then  ${\cal D}(\Gamma)$ is $N$-nicely affine with $r$ parallel classes of blocks, and $\Gamma$ has diameter $4$. 
\end{proposition}
\proof
Assume that $s\ge 4$, so that $\diam(\Gamma)\geq 4$.
Let $\Sigma$ be the graph with $V\Sigma= B'$, such that for any two $w,w' \in B'$, $\{w,w'\}\in E\Sigma$  if and only if $\dist_\Gamma(w,w')=2$ in $\Gamma$. 
Then $G_{ B'}=G$ induces a vertex-transitive group of automorphisms of $\Sigma$ such that
$\Sigma$ is $(G, [s/2])$-distance transitive, by \cite[Lemma 2.7]{DGLP}, and ${\cal N}\setminus \{B\}$ is the set of 
$N$-orbits in $ B'$ forming a nontrivial $G$-invariant partition of $V\Sigma$.

Let  $C\in {\cal N}\setminus \{B\}$ and $x\in C$. Since $\Sigma$ is connected, there exists $x'\in B'$ with $\dist_\Gamma(x,x')=2$.  By Lemma \ref{lem:diffneigh}, $\Gamma(x)\ne\Gamma(x')$, and so there is a 2-arc $(x,v,x')$ in $\Gamma$ 
and a vertex $v'\in\Gamma(x')\setminus\Gamma(x)$.
By Lemma \ref{cor:star}, $\Gamma(v')\cap C=\{y\}$ for some $y$, and we have a 4-arc $(x,v,x',v',y)$ where $\dist_\Gamma(C)\geq 4$. Hence $\dist_\Gamma(x,y)=4$, $\dist_\Sigma(x,y)=2$, and so  $\dist_\Sigma(C)= 2$.
Since $s\geq4$, Lemma 5.2 of \cite{DGLP} applies to $\Sigma$, with $G$-invariant partition ${\cal N}'= {\cal N}\setminus \{B\}=\{C_1,C_2,\dots,C_r\}$.
Since $r\geq3$ 
and  $G$ is transitive on $ B'=V\Sigma$, 
 $\Sigma$ must be complete multipartite with $r$ blocks of size $\ell$ (denoted by  $\Sigma\cong\bfK_{r[\ell]}$) and $\Sigma_{{\cal N}'}\cong\bfK_r$. Thus, for $x,z\in  B'$ lying in distinct blocks of ${\cal N}'$, we have $\dist_\Gamma(x,z)=2$,
and $\dist_\Gamma(x,y)=4$ for distinct vertices $x,y\in  B'$ in the same block of ${\cal N}'$. Therefore $\diam_\Gamma(x):=max\{\dist_\Gamma(x,y)|y\in V\Gamma\}$ is equal to $4$ or $5$, and as $s\geq 4$, $G_x$ is transitive on $\Sigma_i(x)=\Gamma_{2i}(x)$ for $i=1,2$,  for all $x \in  B'$.
In terms of $\mathcal{D}(\Gamma)$ ( which is resolvable by Proposition \ref{star quotient}), this means that  $G$ is transitive om intersecting and non-intersecting block pairs, and so any two `blocks' in different parallel classes meet in a constant number $\mu$ of `points'. Therefore  $\mathcal{D}(\Gamma)$ is an affine resolvable design.
Since $N$ is transitive on points and the $N$-orbits on blocks are exactly the parallel classes, we have that  $\mathcal{D}(\Gamma)$ is $N$-nicely affine. 

Let $C'\in {\cal N}\setminus \{B,C\}$ and $v\in B\setminus\Gamma_1(x)$. Then $\Gamma(v)\cap C'=\{y'\}$,say, by Lemma \ref{cor:star}, and we have shown that $\dist_\Gamma(x,y')=2$.
 Hence $v\in\Gamma_3(x)$, so $\diam_\Gamma(x)=4$.
This also implies that for $v\in B$, we have $B'\subseteq \Gamma_1(v)\cup\Gamma_3(v)$, hence $\diam_\Gamma(v)=3$ or $4$ (note that it is possible for $\Gamma_4(v)$ to be empty, in which case there are no non-collinear point pairs). We conclude that $\diam(\Gamma)=4$, and hence $s=4$. 
\qed

\noindent
{\it Proof of Theorem $\ref{main}$.} 
Suppose that (a) holds. Then $\Gamma$ satisfies Hypothesis \ref{hyp} with $r\geq 3$ and $s=4$. By Proposition \ref{prop:affine},  
${\cal D}(\Gamma)$ is $N$-nicely affine with $r$ parallel classes of blocks, and by Proposition \ref{4dist=PT}, 
$\mathcal{D}(\Gamma)$ is $G$-pairwise transitive.

Now assume that (b) holds. Then in particular $\mathcal{D}(\Gamma)$ is an affine design, and so by Proposition \ref{4dist=PT}, 
$\Gamma$ is locally $(G,4)$-distance transitive. 
Since ${\cal D}(\Gamma)$ is $N$-nicely affine, it follows immediately that  $\Gamma$ is $r$-starlike relative to $N$.
\qed

\noindent
{\it Proof of Theorem $\ref{properties}$.}
 Suppose $\Gamma$ and $G$ satisfy  the equivalent conditions of Theorem \ref{main} for some $N\lhd G$.
Then $\Gamma$ is not complete bipartite, so  $G^B$ and $G^{B'}$  are faithful by Lemma \ref{lem:diffneigh}.
By Lemma \ref{lem:pt}, $G$ has rank 2 or 3 on points and has rank 3 and is imprimitive on blocks. 
Since ${\cal D}(\Gamma)$ is $N$-nicely affine, for any point $x$, $G_x$ is 2-transitive on the set of $N$-orbits on blocks. Since the vertex $x$ is adjacent to exactly one vertex in each $N$-orbot on blocks, the final statement holds. 
\qed

\section{Graphs in Theorem \ref{main}: intersection arrays }

A graph satisfying the equivalent conditions of Theorem \ref{main} is actually locally $G$-distance transitive since by Proposition \ref{prop:affine} it has diameter $4$. 
In particular, it is a distance biregular graph, see \cite{GodST}  and \cite{biregular}. Thus it has  a pair of {\it  intersection arrays}. 
Let $x$ be a vertex of $\Gamma$ and $i\le \diam(\Gamma)$. Let $y\in \Gamma_i(x)$, 
then let $a_i(x):=|\Gamma_1(y)\cap\Gamma_{i}(x)|$, $b_i(x):=|\Gamma_1(y)\cap\Gamma_{i+1}(x)|$ (if $i<\diam_\Gamma(x)$),  and  $c_i(x):=|\Gamma_1(y)\cap\Gamma_{i-1}(x)|$ (if $i>0$). The numbers  $a_i(x)$, $b_i(x)$, and $c_i(x)$ do not depend on the choice of $y$. Notice that $a_i(x)+b_i(x)+c_i(x)$ is equal to the valency of $y$, and so provided this valency is known, $a_i(x)$ can be deduced from the other two numbers. Notice also that $b_0(x)$ is the valency of $x$.
We now define the intersection array of $\Gamma$ at $x$ by 
\[
 \iota(\Gamma,x)=(b_0(x),b_1(x),\dots,b_{\diam(\Gamma)}(x); c_1(x),c_2(x),\dots,c_{\diam(\Gamma)}(x)).
\]

We denote by $\iota(\Gamma)$ the array  $\iota(\Gamma,x)$ where $x\in B$, and by  $\iota'(\Gamma)$ the array  $\iota(\Gamma,x')$ where $x'\in B'$. 
Since $b_i(x)$ only depends on whether $x$ is in $B$ or in $B'$, we write $b_i:=b_i(x)$ for $x\in B$ and $b'_i:=b_i(x')$ for $x'\in B'$. We define $c_i$ and $c'_i$ similarly.

\begin{proposition}\label{arrays}
Assume Hypothesis $\ref{hyp}$ holds and $s\geq 4$. 
Let ${\cal N}:=\{B,C_1,C_2,\ldots C_r\}$,  let  $|C_i|=\ell$, and let $k$ be the valency of the vertices in the bipart $B'$.
Then  $\ell$ divides $k$ and we have the following intersection array: $$\iota'(\Gamma)=(k,r-1,\frac{k(\ell-1)}{\ell},1;1,\frac{k}{\ell},r-1,k).$$
Moreover one of the following happens:
\begin{itemize}
\item $\diam_\Gamma(v)=3$ for $v\in B$, $\ell-1$ divides $k-1$, $r=(k\ell-1)/(\ell-1)$ and the other intersection array is 
$$\iota(\Gamma)=(r,k-1,k;1,r-k,k);$$ 
\item $\diam_\Gamma(v)=4$ for $v\in B$, $\ell(k-1)$ divides $k(r-1)(\ell-1)$, $r(k-\ell)+k(\ell-1)$ divides $k(k-1)(r-1)$, 
and the other intersection array is 
$$\iota(\Gamma)=(r,k-1,\frac{(r-1)k(\ell-1)}{\ell(k-1)},\frac{k(k\ell-\ell r+r-1)}{r(k-\ell)+k(\ell-1)};$$ $$\hspace*{3cm} 1,\frac{r(k-\ell)+k(\ell-1)}{\ell(k-1)},\frac{k(k-1)(r-1)}{r(k-\ell)+k(\ell-1)},r).$$
\end{itemize}
\end{proposition}
\proof
By Proposition \ref{prop:affine}, $\Gamma$ is locally $G$-distance transitive.  Since $\Gamma$ is bipartite, all the $a_i$'s and $a'_i$s are equal to $0$.
In \cite{DGLP} we proved that  $c_{i}c_{i+1}=c'_{i}c'_{i+1}$ for $i$ even (*), and $b_{i}b_{i+1}=b'_{i}b'_{i+1}$ for $i$ odd (**).

We first determine $\iota'(\Gamma)$. Let $x\in B'$, say $x\in C_1$. We have seen that $\diam_\Gamma(x)=4$.
It is an easy consequence of what we have seen above that $|\Gamma_1(x)|=k$, $|\Gamma_2(x)|=(r-1)\ell$, $|\Gamma_3(x)|=k(\ell-1)$, and $|\Gamma_4(x)|=\ell-1$.
We have already seen that $b_0=r$ and $b'_0=k$.
Obviously $c'_1=1$, and so $b'_1=r-1$.
Counting the edges between $\Gamma_1(x)$ and $\Gamma_2(x)$ we have $c'_2=\frac{|\Gamma_1(x)|b'_1}{|\Gamma_2(x)|}=\frac{k}{\ell}$, and therefore  $b'_2=k-\frac{k}{\ell}=\frac{k(\ell-1)}{\ell}$. Hence we get the divisibility condition that $\ell$ divides $k$. Note that $c'_2$ is also equal to $\mu$ (the number of points incident with two non-parallel blocks).
Counting the edges between $\Gamma_2(x)$ and $\Gamma_3(x)$ we have $c'_3=\frac{|\Gamma_2(x)|b'_2}{|\Gamma_3(x)|}=r-1$, and therefore $b'_3=r-(r-1)=1$. Since all vertices adjacent to a point in $\Gamma_4(x)$ are in   $\Gamma_3(x)$, we have $c'_4=k$ and $b'_4=0$.
Hence  $\iota'(\Gamma)=(k,r-1,\frac{k(\ell-1)}{\ell},1;1,\frac{k}{\ell},r-1,k)$.

Now we determine  $\iota(\Gamma)$. Let $v\in B$. We have seen that $\diam_\Gamma(v)=3$ or $4$.
Obviously $c_1=1$, and so $b_1=k-1$.
By (**), we have $b_2=\frac{b'_1b'_2}{b_1}=\frac{(r-1)k(\ell-1)}{\ell(k-1)}$.  Hence we get the divisibility condition that $\ell(k-1)$ divides $k(r-1)(\ell-1)$.
Therefore, $c_2=r-b_2=\frac{r(k-\ell)+k(\ell-1)}{\ell(k-1)}$.
By (*), we have $c_3=\frac{c'_2c'_3}{c_2}=\frac{k(k-1)(r-1)}{r(k-\ell)+k(\ell-1)}$ which adds yet another divisibility condition. Therefore, $b_3=k-c_3=\frac{k(k\ell-r\ell+r-1)}{r(k-\ell)+k(k-1)}$.

Assume that $b_3$ is equal to $0$, that is $\diam_\Gamma(v)=3$. In that case $b_4=c_4=0$ and $k\ell-r\ell+r-1=0$. Hence $r=\frac{k\ell-1}{\ell-1}$ and  $\ell-1$ divides $k-1$, which imply the two divisibility conditions of the preceding paragraph.
 Moreover, we can compute the simpler expressions  $b_2=k$, $c_2=\frac{k-1}{\ell-1}=r-k$ and $c_3=k$. 

Assume now that $b_3\neq 0$. Then $\diam_\Gamma(v)=4$, 
$b_{4}=b'_{3}b'_{4}/b_{3}=0$ by (**), and so  $c_{4}=b'_0-b_{4}=r$.

Therefore we get $\iota(\Gamma)$ as stated.
\qed

\section{Regular point action by $N$}

In this section, we will prove Theorem \ref{thm:regN} through a series of lemmas.
In all the following lemmas, we assume that ${\cal D}$ is  a $G$-pairwise transitive and $N$-nicely affine design 
with $r$ parallel classes, for some $r\geq 3$ and  $N\lhd G$ (so that $\Gamma=\Gamma(\calD)$ and $\calD$ satisfy Theorem \ref{main}(a) and (b)).
In addition we assume that $N$ acts regularly the point set, that is to say, $N$ is transitive and only the identity fixes a point.  

Let $C_1,\ldots,C_r$ be the orbits of $N$ on blocks. Then the blocks in each $C_i$ partition the point set of $\mathcal{D}$, 
giving a resolution of the design. Since  $N$ is regular and faithful on the point set, by Lemma~\ref{lem:diffneigh}, 
we may identify the point set with $N$ in such a way that $N$ acts by right multiplication. 
Let $v$ be the point identified with $1_N$. 
Then $G=N.G_v$, where $G_v$ acts on $N$ by conjugation (see \cite[p.9]{Cam99}).
By Theorem \ref{properties}, we know that  $G$ has rank at most 3 in its action on $N$.

\begin{lemma}\label{lem:subgpsN} There are $r$ subgroups $M_1,\ldots, M_r$ of $N$ such that the blocks of the designs  $\mathcal{D}$ can be identified with $\{M_ix|i=1,\ldots, r, x\in N\}$. Moreover $N$ acts by right multiplication and $G_v$ by conjugation on thoses cosets. In particular, $G_v$ leaves $\{M_1,\ldots, M_r\}$ invariant.
\end{lemma}
\proof 
Each  block of the design  can be identified with its point-set, that is, a set of elements of $N$. Since $G$ is transitive on pairs of intersecting blocks,
no two blocks have the same point-set, and so this identification is not ambiguous.
Let $w_i$ be the block of $C_i$ containing $1_N$. Let $M_i\subset N$ be its point-set. We claim that $M_i$ is a subgroup of $N$. 
By definition, $M_i$ contains $1_N$. Take $m\in M_i$. Then  $w_i^m$ is a block disjoint from $M_i$ or $M_i$ itself. Since $N$ acts on itself by right multiplication, $w_i^m$ is the
subset $M_im$ of $N$. Now $m\in M_i\cap M_im$, and so $M_im=M_i$. This implies that $M_i$ is closed under multiplication. Also  $1_N\in M_im^{-1}\cap M_i$, so $M_im^{-1}= M_i$. 
This implies that $M_i$ is closed  under taking inverses. Thus $M_i$ is a subgroup of $N$, for each $i$. Now $C_i$ is an $N$-orbit, so for $z_i\in C_i$, there exists $n\in N$ such that 
$z_i=w_i^n$. Using the identification we get that the point-set of $z_i$ is $M_in$ (recall that $N$ acts by right multiplication on $N$), so $z_i$ is identified with a right coset of $M_i$ in $N$. It follows that $C_i$ can be identified with the set of right cosets of $M_i$ in $N$. 

It follows from the identification made that $N$ acts by right multiplication. 
Let $h\in G_v$. We know $w_i^h=w_j$ for some $j=1,\ldots,r$. Since $M_i$ is the point-set of $w_i$, we have that $M_j$ is the point-set of $w_j=w_i^h$, so is equal to $M_i^h$. Hence $G_v$ permutes $\{M_1,\ldots, M_r\}$, acting by conjugation.
Then $h$ maps the coset $M_in=\{mn|m\in M_i,n\in N\}$ to $\{m^hx^h|m\in M_i,n\in N\}=M_i^hx^h=M_jx^h$ which is also a coset but of $M_i^h$. Thus $G_v$ acts by conjugation on the cosets.
\qed

\begin{lemma}\label{lem:pgroups} We have that $G_v$ is $2$-transitive on $\{M_1,\ldots, M_r\}$, where $v=1_N$. 
Moreover, the subgroups $M_i$ are isomorphic $p$-groups of exponent $p$.
\end{lemma}
\proof
Let $v=1_N$.
By  Theorem \ref{properties}, $G_v$ is $2$-transitive on $\Gamma_1(v)$, and as explained in the previous proof each element of $\Gamma_1(v)$ is unambguously associated to one subgroup $M_i$, so $G_v$ is $2$-transitive on $\{M_1,\ldots, M_r\}$. Recall that each $N$-orbit is the set of cosets of an $M_i$. We choose in each orbit the representant $M_i$.   Note that an element in $G_v$ necessarily maps $M_i$ to some $M_j$, so  $G_v$ is $2$-transitive on $\{M_1,\ldots, M_r\}$.

Notice in particular that all the $M_i$ are isomorphic. Since ${\cal D}$ is transitive on pairs of collinear points, we also have that $G_v$ acts transitively by conjugation on $\cup_{i=1}^rM_i\setminus\{1_N\}$.
Thus all non-identity elements of all the $M_i$ must have the same prime order $p$. Hence the $M_i$ are $p$-groups of exponent $p$.
\qed

\begin{lemma}\label{lem:Npgroup} We have $N=M_1M_2$ is a $p$-group satisfying $|N|\leq |M_1|^2$. 
\end{lemma}
\proof
Since ${\cal D}$ is an affine design, we have that $M_2$ meets all cosets of $M_1$. Since the cosets of $M_1$ partition $N$, it follows that every element of $N$ is in $M_1m_2$ for some $m_2\in M_2$, and so $N=M_1M_2$. Thus $|N|=\frac{|M_1||M_2|}{|M_1\cap M_2|}$ is a $p$-power, since $M_1$ and $M_2$ are both $p$-groups by Lemma \ref{lem:pgroups}. Moreover $|N|\leq |M_1||M_2|=|M_1|^2$.
\qed

{\it Proof of Theorem \ref{thm:regN}.}
Suppose $N$ is not abelian, then $N$ has at least one proper non-trivial characteristic subgroup, namely $Z(N)$ (which is not trivial since $N$ is a $p$-group). Let $A$ be any proper non-trivial characteristic subgroup of $N$. Since $G_v\leq \Aut(N)$, it follows that $G_v$ stabilises $A$ setwise.
 Then $G$ has rank $3$ on $N$ and one of $A\setminus\{1_N\}$ and $N\setminus A$ is the set of points collinear with $v$ and the other is the set of points not collinear with $v$. Suppose  $A\setminus\{1_N\}$  is the set of points collinear with $v$. Then $A=\cup_{i=1}^rM_i$. 
Since $N=M_1M_2$, this implies that $N\subseteq A$, a contradiction. Thus $A\setminus\{1_N\}=N\setminus \cup_{i=1}^rM_i$ and $N\setminus A=\cup_{i=1}^rM_i\setminus\{1_N\}$. 
Therefore $A$ is uniquely determined and $G$ has only one proper non-trivial characteristic subgroup, that is, $Z(N)=\{1_N\}\cup (N\setminus \cup_{i=1}^rM_i)$.
It follows that $M_i\cap Z(N)=\{1_N\}$ for all $i$. This means that the natural homomorphism $\nu$ from $N$ to $N/Z(N)$ restricted to $M_i$ is injective, and so $M_i\cong \nu(M_i)$.
The second center $Z_2(N)$ is also characteristic and contains $Z(N)$ strictly (since $N/Z(N)$ is also a $p$-group). Thus we must have $Z_2(N)=G$, and so $N/Z(N)$ is abelian. It follows that  $M_i$ is abelian for all $i$ since $\nu(M_i)$ is. 
Then $N=M_1M_2$ is a product of two abelian groups, and by a result of Ito \cite{Ito}, either $M_1$ or $M_2$ contains 
a nontrivial normal subgroup $S$ of $N$. But then $S$ meets $Z(N)$ nontrivially, and so $M_1$ or $M_2$ meets the center nontrivially, a contradiction. 
We conclude that $N$ must be abelian.

Since $M_1$ and $M_2$ have exponent $p$ by Lemma \ref{lem:pgroups}, and $N=M_1M_2$ is abelian, it follows that $N$ itself has exponent $p$ and so $N$ is elementary abelian.
Therefore we are now going to call its neutral element $0$ and write $N$ additively. Note that if $|N|=p^d$ then $N$ can be identified with a vector space $V=V(d,p)$, subgroups of $N$ are just subpaces, and $G\leq \AGL(d,p)$.
Since the points of the design are identified with $N$, we have that the point set of ${\cal D}$ can be identified with the vector space $V=V(d,p)$ and the blocks can be identified with cosets of subspaces, by Lemma \ref{lem:subgpsN}. 

Now we consider the conditions of Construction \ref{constr-regN}.
We  proved (a) in Theorem \ref{properties}, (b) in Lemmas \ref{lem:subgpsN} and \ref{lem:pgroups}, and (c) in Lemma \ref{lem:Npgroup}.  
We now prove (d). Let $x+M_1$ and $y+M_1$ be two non-trivial cosets of $M_1$ (so $x,y\notin M_1$). Then $x+M_1,y+M_1$ do not contain the point $0$, and by Lemma \ref{lem:subgpsN}, $G_0$ acts by conjugation on cosets. Since ${\cal D}$ is transitive on non-incident point-block pairs, there exists $g\in G_0$ such that $y+M_1=x^g+M_1^g$. Hence $M_1^g=M_1$ and (d) follows. 
If $\cup_{i=1}^rM_i=V$, then the set of pairs of non-collinear points is empty and so $G$ is 2-transitive on $V=N$ since $G_0$ is transitive on the set of all pairs of points. Hence (e) holds.
\qed

\section{The case $r=2$}
In this section, we prove Theorem \ref{r=2}.
 Suppose throughout that $\Gamma$ is not complete bipartite, and is locally $(G,s)$-distance transitive
and is 2-starlike relative to a normal subgroup $N$ of $G$, for some $s$ such that $2\leq s\leq \diam(\Gamma)$, that is to say, Hypothesis \ref{hyp} holds with $r=2$. We use that notation.
Recall that the subdivision graph $S(\Sigma)$ of a graph $\Sigma$ 
is the bipartite graph with ordered bipartition $(E\Sigma|V\Sigma)$ and adjacency given by containment.   

\begin{lemma}\label{subdgraph}
 
The graph  $\Gamma$ is the subdivision graph of a connected, bipartite graph $\Sigma$ such that $\diam(\Gamma)=2\,\diam(\Sigma)$. 
\end{lemma}

\proof
Let ${\cal N}=\{B,C_1,C_2\}$, as in Hypothesis~\ref{hyp}. By Lemma \ref{cor:star}, $G=G^+$.
Each vertex in $B$ is adjacent to exactly two vertices in $B'$, by Lemma \ref{cor:star}(c).
By Lemma \ref{lem:diffneigh}, $\Gamma(x)\ne\Gamma(y)$ for distinct $x,y\in B$. 
Therefore two vertices of $ B'$  at distance 2 have a unique common neighbour.
Let $\Sigma$ be the graph  with $V\Sigma= B'$ and edges the pairs from $B'$ at distance 2 in $\Gamma$. 
Then each edge of $\Sigma$ is of the form $\Gamma(x)$ for a unique $x\in B$, and hence $\Gamma$ is the subdivision graph of $\Sigma$. 
Since $\Gamma$ is connected, so is $\Sigma$.
Moreover, $ B'$ is split into two orbits $C_1$, $C_2$ by $N$ and each element of $B$ is adjacent to one vertex in each orbit. 
In other words, the graph $\Sigma$ is bipartite. Finally, $\Sigma$ is regular since $G$ is transitive on $B'$, and 
the assertion about diameters of these graphs follows from \cite[Lemma 5]{DD}.
\qed

\begin{lemma}\label{cycles}
   If the graph $\Sigma$ in Lemma~$\ref{subdgraph}$ has valency $2$, then $\Sigma=C_{2\ell}$ and $\Gamma=C_{4\ell}$ for some $\ell\geq s/2$.
\end{lemma}

\proof
Since $\Sigma$ is connected and bipartite of valency 2 it follows that $\Sigma=C_{2\ell}$ for 
some $\ell$. The subdivision graph of $\Sigma$ is a cycle of twice the length, and hence by
Lemma~\ref{subdgraph}, $\Gamma=C_{4\ell}$. Finally $s\leq\diam(\Gamma)=2\ell$ by definition of local $s$-distance transitivity.
\qed

\noindent
{\it Proof of Theorem $\ref{r=2}$} 
Suppose first that $\Gamma$ is the complete bipartite graph $\bfK_{n,m}$ with $N$ transitive on the part $B$ of size $n$ and having two orbits $C_1, C_2$ on the part of size $m$. Then $2=\diam(\Gamma)$, so $s=2$. Also, since $G$ is transitive on $C_1\cup C_2$ it follows that $|C_1|=|C_2|$ and so $m$ is even. Moreover, $N\leqslant S_n\times (S_{m/2}\times S_{m/2})$ and $G\leqslant N_{\Aut(\Gamma)}(N)\leqslant S_n\times (S_{m/2}\wr S_2)$. Given $v\in C_1$ the elements at distance 2 from $v$ are the elements of $C_2\cup (C_1\backslash\{v\})$. Thus $S_n\times (S_{m/2}\wr S_2)$ is not locally 2-distance transitive when $m\geq 4$  and so $m=2$.

From here on we suppose $\Gamma$ is not complete bipartite and that  $\Gamma, G, s$ satisfy the hypotheses of Theorem~\ref{r=2}, 
so $\Gamma$ satisfies 
the assumptions above for this section, and Hypothesis \ref{hyp} holds with $r=2$. By Lemma~\ref{subdgraph}, 
$\Gamma$ is the subdivision graph of a connected, regular bipartite graph $\Sigma$ with $2\diam(\Sigma)=\diam(\Gamma)$, 
and we may assume that $B=E\Sigma$, $B'=V\Sigma$.
If $\Sigma$ has valency 2, then  $\Gamma=C_{4\ell}$ for some $\ell\geq s/2$, by Lemma \ref{cycles}, and there is nothing further to be proved. 
So assume that $\Sigma$ has valency $k\geq3$. If $s<\diam(\Gamma)=2\diam(\Sigma)$ then, by \cite[Theorem 1.2]{DDP12}, 
$\Sigma$ is  $(G, \lceil \frac{s+1}{2}\rceil)$-arc transitive. In this case, by Weiss' Theorem  \cite{8-arc-trans}, 
we have that $\lceil\frac{s+1}{2}\rceil\leq 7$, and so $s\leq 13$.

Finally suppose that $s=\diam(\Gamma)=2\diam(\Sigma)$. Then it follows that $\Sigma$ is one 
of the graphs in Remark~\ref{r=2list}(c) from \cite[Theorem 1 and Corollary 2]{DD}
(we simply remove the non-bipartite graphs from the lists for these results.)
Since $\Gamma$ is in particular
$(G,2\diam(\Sigma)-1)$-distance transitive, it follows from 
\cite[Theorem 1.2]{DDP12} that $\Sigma$ is $(G,\diam(\Sigma))$-arc transitive. Applying Weiss's Theorem again, 
$\diam(\Sigma)\leq 7$, whence $s=2\diam(\Sigma)\leq 14$.
\qed

\begin{lemma}\label{lem:last}
If  $\Gamma=S(\Sigma)$ is as in Theorem $\ref{r=2}(ii)$, then:

(a) if $\Sigma$ is  $(G, \lceil \frac{s+1}{2}\rceil)$-arc transitive, then $\Gamma$ is locally $(G,s)$-distance transitive;

(b) if $\Sigma$ is in the list of Remark $\ref{r=2list}(c)$, then there exists an automorphism group $G$ such that $\Gamma$ is locally $(G,\diam(\Gamma))$-distance transitive;

(c) $\Gamma$ is $2$-starlike relative to a normal subgroup $N$ of $G$, where $G$ is as in part (a) or (b).
\end{lemma}
\proof
(a) Suppose $\Sigma$ is  $(G, \lceil \frac{s+1}{2}\rceil)$-arc transitive, then  $\Gamma$ is locally $(G,s)$-distance transitive by Theorem 1.2 in \cite{DDP12}.

(b)  Suppose $\Sigma$ is in the list of Remark \ref{r=2list}(c). For each of these graphs, suitable groups $G$ are specified in \cite{DD}, for which $\Gamma$ is locally $(G,\diam(\Gamma))$-distance transitive. 

(c) Let $G$ be as in part (a) or (b), and let $N$ be the subgroup of $G$ stabilising the biparts of $\Sigma$ setwise. 
In all cases, $G$ is transitive on the set of vertices of the bipartite graph $\Sigma$, so $N$ has two orbits on $V\Sigma=B'$.
We claim that $N$ is transitive on $E\Sigma$, and hence transitive on $B$.
Let $e=\{u,v\}$ and $e'=\{u',v'\}$ be two edges of $\Sigma$. Since $N$ is transitive on each bipart of $\Sigma$, we can assume that $u=u'$. Then both $e$ and $e'$ are vertices at distance 1 from $u$ in $\Gamma$, and since $\Gamma$ is locally $(G,1)$-distance transitive, there is an element of $G$ fixing $u$ and mapping $e$ to $e'$. This element is in $N$ since it preserves the bipart containing $u$, and so the claim is proved. 
Hence $\Gamma$ is 2-starlike relative to $N$. 
\qed

\section*{Acknowledments}
The authors wish to thank the reviewers for helpful suggestions on the exposition.

\end{document}